\documentclass{article}[12pt]
\usepackage{epsfig,makeidx,amsmath,amsfonts,latexsym,subfigure}
\newtheorem{theorem}{Theorem}

\newtheorem{lemma}{Lemma}

\begin{document}

\title{\large\textbf{{The initial configuration is irrelevant for
the possibility\linebreak of mutual unbounded growth in
the\linebreak two-type Richardson model}}}

\author{Maria Deijfen \thanks{Stockholm University. E-mail: mia@matematik.su.se}
\and Olle H\"{a}ggstr\"{o}m \thanks{Chalmers University of
Technology. E-mail: olleh@math.chalmers.se}}

\date{May 2003}

\maketitle

\thispagestyle{empty}

\begin{abstract}

\noindent The two-type Richardson model describes the growth of
two competing infections on $\mathbb{Z}^d$. At time 0 two disjoint
finite sets $\xi_1,\xi_2\subset \mathbb{Z}^d$ are infected with
type 1 and type 2 infection respectively. An uninfected site then
becomes type 1 (2) infected at a rate proportional to the number
of type 1 (2) infected nearest neighbors and once infected it
remains so forever. The main result in this paper is, loosely
speaking, that the choice of the initial sets $\xi_1$ and $\xi_2$
is irrelevant in deciding whether the event of mutual unbounded
growth for the two infection types has positive probability or
not.

\vspace{1cm}


\vspace{0.5cm}

\hspace*{4.8cm} Secondary 82B43.
\end{abstract}

\section{Introduction}

\noindent This paper is concerned with certain models for random
growth and competition on the cubic lattice $\mathbb{Z}^d$ in
dimension $d\geq 2$. The Richardson model, introduced in
Richardson (1973), is one of the simplest models for growth on
$\mathbb{Z}^d$. Each site is in either of two states, denoted by 0
and 1, and the set of sites in state 1 increases to cover all of
$\mathbb{Z}^d$. The dynamics is that a site in state 0 is
transferred to state 1 at a rate proportional to the number of
nearest neighbors in state 1 and a site in state 1 remains there
forever. This is equivalent to first-passage percolation with
i.i.d.\ exponential passage times.\medskip

\noindent In H\"{a}ggstr\"{o}m and Pemantle (1998), a generalized
version of the Richardson model is introduced. There are now three
possible states, denoted by 0, 1 and 2 respectively, for each
site. A site in state 0 is transferred to state 1 (2) with rate
$\lambda_1$ ($\lambda_2$) times the number of nearest neighbors in
state 1 (2) and the states 1 and 2 are absorbing. Here
$\lambda_1,\lambda_2>0$ are the parameters of the model. For
disjoint sets $\xi_1,\xi_2\subset \mathbb{Z}^d$, let
$P^{\lambda_1,\lambda_2}_{\xi_1,\xi_2}$ denote the probability law
of the generalized process started at time zero with the sites in
$\xi_1$ being in state $1$, the sites in $\xi_2$ in state 2 and
the rest of $\mathbb{Z}^d$ in state 0. The states 1 and 2 may be
thought of as representing two different types of infection and
the model then describes the growth of two infections competing
for space on $\mathbb{Z}^d$. There are two possible scenarios:
Either one of the infection types at some point completely
surrounds the other, preventing the surrounded type from growing
any further, or both infection types keep growing indefinitely.
Write $A$ for the latter event, that is,
$$
A=\{\textrm{both infection types reach sites arbitrarily far away
from the origin}\}.
$$
\noindent If the initial sets $\xi_1$ and $\xi_2$ are finite,
clearly $A^c$ has positive probability. To decide whether $A$ has
positive probability or not is more intricate. First, assume that
$\xi_1=\{\textbf{0}\}$ and $\xi_2=\{\textbf{1}\}$, where
$\textbf{n}=(n,0,\ldots,0)$. Intuitively, $A$ should in this case
occur with positive probability if and only if
$\lambda_1=\lambda_2$. This intuition is partly confirmed in
H\"{a}ggstr\"{o}m and Pemantle (1998,2000). The main result in the
first paper is that, if $\lambda_1=\lambda_2$ and $d=2$, then
$P^{\lambda_1,\lambda_2}_{\textbf{0,1}}(A)>0$. In the second paper
it is proved that, if $d\geq 2$ and $\lambda_1$ is held fixed,
then $P^{\lambda_1,\lambda_2}_{\textbf{0,1}}(A)=0$ for all but at
most countably many values of $\lambda_2$. The aim in the present
paper is to show that the choice of initial sets is basically
irrelevant for these results. Of course, if one set completely
surrounds the other, then mutual unbounded growth is ruled out. To
formulate our main result, we therefore employ the following
definition.\medskip

\noindent \textbf{Definition 1} Let $\xi_1$ and $\xi_2$ be two
disjoint finite subsets of $\mathbb{Z}^d$. We say that one of the
sets ($\xi_i$) \emph{strangles} the other ($\xi_j$) if there
exists no infinite self-avoiding path in $\mathbb{Z}^d$ that
starts at a vertex in $\xi_j$ and that does not intersect $\xi_i$.
The pair $(\xi_1,\xi_2)$ is said to be \emph{fertile} if neither
of the sets strangles the other.\medskip

\noindent The main result is as follows.

\begin{theorem}\label{th:huvudresIII}
Let $(\xi_1,\xi_2)$ and $(\xi_1',\xi_2')$ be two fertile pairs of
disjoint finite subsets of $\mathbb{Z}^d$. For all choices of
$(\lambda_1,\lambda_2)$, we have
$$
P^{\lambda_1,\lambda_2}_{\xi_1,\xi_2}(A)>0\Leftrightarrow
P^{\lambda_1,\lambda_2}_{\xi_1',\xi_2'}(A)>0.
$$
\end{theorem}
\noindent This implies that mutual unbounded growth has positive
probability when starting from $(\{\textbf{0}\},\{\textbf{1}\})$
if and only if it occurs with positive probability for every other
fertile initial configuration as well. Hence the results in
H\"{a}ggstr\"{o}m and Pemantle (1998,2000) extend to arbitrary
initial sets, as desired.\medskip

\noindent H\"{a}ggstr\"{o}m and Pemantle (1998) contains a special
case of Theorem \ref{th:huvudresIII}, namely the case when $d=2$,
$\lambda_1=\lambda_2$ and $(\xi_1,\xi_2)$ and $(\xi_1',\xi_2')$
consist of single sites. The proof readily extends to the case
where $\lambda_1$ and $\lambda_2$ are arbitrary and $\xi_1$ and
$\xi_2$ are both connected sets. However, the proof fails to
extend to more general initial configurations, and, since it uses
planarity, it is unclear whether it extends to the case $d\geq 3$.
These difficulties are overcome in the present paper.\medskip

\noindent We mention also that a related model for competing
growth in continuous space was studied by Deijfen \emph{et al}
(2003), and the results obtained there include a kind of continuum
analogue of our main result.

\section{Preliminaries}

\noindent In this section we give a concrete construction of the
two-type Richardson model that suits our purposes. We also
introduce some notation and formulate a lemma that will be
important in the proof of Theorem \ref{th:huvudresIII}.\medskip

\noindent To begin with, note that, by time-scaling and symmetry,
we may restrict our attention to two-type processes with rates
$(1,\lambda)$ for some $\lambda\leq 1$. To build up such a
process, define the distance between two sites
$x=(x_1,\ldots,x_d)$ and $y=(y_1,\ldots,y_d)$ on $\mathbb{Z}^d$ by
$\delta(x,y)=\sum_{i=1}^d|x_i-y_i|$, and call two sites nearest
neighbors if they are located at distance 1 from each other.
Independently for each ordered pair $(x,y)$ of nearest neighbor
sites on $\mathbb{Z}^d$, associate a unit rate Poisson process
$P^{(x,y)}$ and, for $\lambda\in[0,1]$, write $\lambda P^{(x,y)}$
for the thinning of $P^{(x,y)}$ obtained by removing each Poisson
occurrence with probability $1-\lambda$. Intuitively, at the times
of the occurrences in the process $P^{(x,y)}$, we imagine that a
channel between $x$ and $y$ is opened so that type 1 infection can
be transferred from $x$ to $y$, that is, if $x$ is type 1 infected
at such a time, then $y$ will become type 1 infected as well. The
type 2 infection is controlled analogously by the process $\lambda
P^{(x,y)}$.\medskip

\noindent To formally define the growth process, let
$\Gamma_{n}^i$ denote the set of type $i$ infected sites after $n$
infections and let $T_n$ denote the time point for the $n$th
infection. Also, for a set $\eta\subset\mathbb{Z}^d$, define
$\partial\eta$ to be the set of sites in $\eta$ that has at least
one nearest neighbor in $\eta^c$, that is,
$$
\partial\eta=\{x\in\eta;\hspace{0.1cm}\exists
y\in\mathbb{Z}^d\backslash\eta \textrm{ with } \delta(x,y)=1\}.
$$
The sequences $\{\Gamma_{n}^1\}$, $\{\Gamma_{n}^2\}$ and $\{T_n\}$
are obtained inductively as follows:

\begin{itemize}
\item[1.] Let $\Gamma_{0}^1=\xi_1$, $\Gamma_{0}^2=\xi_2$ and
$T_0=0$.
\item[2.] Given $\Gamma_{n}^1$, $\Gamma_{n}^2$ and $T_n$, define
$T_{n+1}=\min\{\acute{T}_{n+1}^1,\acute{T}_{n+1}^2\}$,
where\smallskip

\noindent $\acute{T}_{n+1}^1=\inf\{T>T_n;\hspace{0.1cm} T\in
P^{(x,y)}$ for some pair $(x,y)$ such that
$x\in\partial\Gamma_{n}^1$
\newline \hspace*{1.7cm} and
$y\not\in\Gamma_{n}^1\cup\Gamma_{n}^2\}$\medskip

\noindent and $\acute{T}_{n+1}^2$ is defined analogously but with
$P^{(x,y)}$ replaced by $\lambda P^{(x,y)}$ and
$\partial\Gamma_n^1$ replaced by $\partial\Gamma_n^2$.

\item[3.] If $T_{n+1}=\acute{T}_{n+1}^1$, then $\Gamma_{n+1}^2=
\Gamma_{n}^2$ and $\Gamma_{n+1}^1=\Gamma_{n}^1\cup\{y\}$, where
$y$ is the site such that $\acute{T}_{n+1}\in P^{(x,y)}$. If
$T_{n+1}=\acute{T}_{n+1}^2$, then $\Gamma_{n}^2$ is updated in the
same way while $\Gamma_{n}^1$ is left unchanged.
\end{itemize}

\noindent The set of type $i$ infected sites at time $t\in
[T_n,T_{n+1})$ is given by $\Gamma_i(t)=\Gamma_{n}^i$ and the
total set of infected sites at time $t$ is
$\Gamma(t)=\Gamma_1(t)\cup \Gamma_2(t)$. When the initial sets
need to be clear from the notation they will be included as
superscripts, for instance $\Gamma^{\xi_1,\xi_2}(t)$ denotes the
set of infected sites at time $t$ in a process started from the
sets $(\xi_1,\xi_2)$. By standard properties of the Poisson
process, the time until an infected site infects an uninfected
nearest neighbor is exponentially distributed and hence
$\{\Gamma(t)\}$ is a Markov process.\medskip

\noindent  We remark that in the original construction of the
two-type Richardson model, given in H\"{a}ggstr\"{o}m and Pemantle
(1998,2000), the type 1 and the type 2 infections are generated by
independent Poisson processes. More precisely, to each ordered
nearest neighbor pair $(x,y)$ two independent Poisson processes,
$P^{(x,y)}_1$ and $P^{(x,y)}_2$, with rates 1 and $\lambda$
respectively, are attached. The process $P^{(x,y)}_1$ then
controls the type 1 infection and $P^{(x,y)}_2$ controls the type
2 infection. However, in the present article we will always assume
that the type 1 and the type 2 infections are generated by the
same Poisson process as described above. Obviously this gives a
growth process with the same distribution as the original
one.\medskip

\noindent One way of describing the evolution of the infection is
to study the \emph{infection graph}, denoted by $\Psi$. It
consists of two disjoint graphs, $\Psi_1$ and $\Psi_2$, describing
the type 1 and the type 2 infection respectively. These are
generated as follows: For $i=1,2$, let $\Psi_i(t)$ be the graph
with vertex set $\Gamma_i(t)$ and edge set obtained by putting an
edge between two nearest neighbor sites $x,y\in\Gamma_i(t)$ if and
only if, at some time $t'\leq t$, $x$ was infected by $y$ or vice
versa. Define $\Psi_i=\lim_{t\rightarrow\infty}\Psi_i(t)$, where
the limit exists since both the vertex set and the edge set of
$\Psi_i(t)$ is increasing in $t$, and let $\Psi=\Psi_1\cup
\Psi_2$. In the one-type Richardson model started from a single
infected site, $\Psi$ is a tree and its features have been studied
in Newman (1995) and H\"{a}ggstr\"{o}m and Pemantle (1998). In
general, $\Psi$ is a forest, that is, each connected component is
a tree.\medskip

\noindent Clearly mutual unbounded growth for the two infection
types in the two-type model occurs if and only if both $\Psi_1$
and $\Psi_2$ contain an infinite path. The following lemma relates
the existence of such paths to the boundary configuration of the
initial set. To formulate it, extend the notation for the
infection graphs to incorporate the initial sets, so that
$\Psi^{\zeta_1,\zeta_2}_i$ denotes the type $i$ infection graph
for a process started from $(\zeta_1,\zeta_2)$.

\begin{lemma}\label{lemma1}
Consider two growth processes with the same infection rates
$(1,\lambda)$, $\lambda\leq 1$, and generated by the same Poisson
processes, but started from two different finite initial
configurations $(\zeta_1,\zeta_2)$ and $(\zeta_1',\zeta_2')$.
Assume that $\zeta_1\cup \zeta_2=\zeta_1'\cup \zeta_2'=\zeta$ and
$\zeta_2\cap\partial\zeta\subset \zeta_2'\cap\partial\zeta$.

\begin{itemize}
\item[\rm{(a)}] If $\Psi^{\zeta_1,\zeta_2}_2$ contains
an infinite path, then so does $\Psi^{\zeta_1',\zeta_2'}_2$.

\item[\rm{(b)}] If there is an infinite path in $\Psi^{\zeta_1,\zeta_2}_1$ starting
at some site $x\in\partial\zeta\cap\zeta_1\cap\zeta_1'$, then the
same path is present in $\Psi^{\zeta_1',\zeta_2'}_1$ as well.
\end{itemize}
\end{lemma}

\noindent \emph{Proof:} Write $\zeta^\circ$ for the interior of
$\zeta$, that is, $\zeta^\circ=\zeta\backslash\partial\zeta$. We
will show that

\begin{equation}\label{eq:1}
\Gamma_1^{\zeta_1,\zeta_2}(t)\backslash \zeta^\circ\supset
\Gamma_1^{\zeta_1',\zeta_2'}(t)\backslash \zeta^\circ
\end{equation}

\noindent and

\begin{equation}\label{eq:2}
\Gamma_2^{\zeta_1,\zeta_2}(t)\backslash \zeta^\circ\subset
\Gamma_2^{\zeta_1',\zeta_2'}(t)\backslash \zeta^\circ
\end{equation}

\noindent for all $t$. To this end, order the time points for the
infections in the two growth processes in one single sequence
$\{\tilde{T}_n\}_{n\geq 0}$, where $\tilde{T}_0:=0$. Note that,
since the growth processes are generated by the same Poisson
processes, infections can take place simultaneously in both
processes. Hence an infection time $\tilde{T}_n$ can represent an
infection that occurs in both processes. Assume that (\ref{eq:1})
and (\ref{eq:2}) hold for $t=\tilde{T}_n$. The only way for
(\ref{eq:1}) to fail at $t=\tilde{T}_{n+1}$ is then that a site
that is uninfected at time $\tilde{T}_n$ in the process started
from $(\zeta_1,\zeta_2)$ becomes type 1 infected in the process
started from $(\zeta_1',\zeta_2')$. However, it is easily seen
that if this should be the case, then the same infection must take
place in the process started from $(\zeta_1,\zeta_2)$ as well.
Hence, both sets $\Gamma_1^{\zeta_1,\zeta_2}
(\tilde{T}_n)\backslash \zeta^\circ$ and
$\Gamma_1^{\zeta_1',\zeta_2'}(\tilde{T}_n)\backslash \zeta^\circ$
are extended by the same site and thus (\ref{eq:1}) is preserved
at time $\tilde{T}_{n+1}$. Analogously it can be seen that
(\ref{eq:2}) is preserved at $\tilde{T}_{n+1}$. Furthermore, by
assumption, the type 2 infected part of $\partial(\zeta_1\cup
\zeta_2)$ is a subset of the type 2 infected part of
$\partial(\zeta_1'\cup \zeta_2')$, implying that (\ref{eq:1}) and
(\ref{eq:2}) hold for $t=\tilde{T}_0$. It follows by induction
over $n$ that the inclusions (\ref{eq:1}) and (\ref{eq:2}) hold
for all $t\in\{\tilde{T}_n\}$ and clearly they must then hold for
all $t\geq 0$. Part (a) follows immediately from
(\ref{eq:2}).\medskip

\noindent To establish (b), we first show that

\begin{equation}\label{eq:3}
\Gamma^{\zeta_1,\zeta_2}(t)\supset
\Gamma^{\zeta_1',\zeta_2'}(t)\textrm{ for all }t\geq 0.
\end{equation}

\noindent To this end, assume that the site $y\not\in \zeta$ is
infected at time $t>0$ in the process started from
$(\zeta_1',\zeta_2')$. If $y$ is type 1 infected, it follows from
(\ref{eq:1}) that $y$ is (type 1) infected at time $t$ in the
process started from $(\zeta_1,\zeta_2)$ as well. So suppose $y$
is type 2 infected. Then there is an infection chain with passage
time at most $t$ that leads from some site in $\zeta_2'$ to $y$.
Clearly, unless some site in the chain becomes type 1 infected
before it is reached by the type 2 infection, the same chain is
present also in the process started from $(\zeta_1,\zeta_2)$.
However, the fact that a site in the chain is type 1 infected can
only decrease the time it takes for the infection to reach $y$.
Hence $y$ is infected at the latest at time $t$ in the process
started from $(\zeta_1,\zeta_2)$ and (\ref{eq:3}) follows.\medskip

\noindent Now let $\{v_n\}$ and $\{e_n\}$ denote the vertex and
edge set respectively of the infinite path starting at $x=v_0$ in
$\Psi^{\zeta_1,\zeta_2}_1$. We will call a site $v_n$
\emph{successful} in the process started from
$(\zeta_1',\zeta_2')$ if it is type 1 infected and infects the
site $v_{n+1}$ via the edge $e_n$ at the latest at the time when
$v_{n+1}$ is infected in the process started from
$(\zeta_1,\zeta_2)$. Using (\ref{eq:3}), it follows easily by
induction over $n$ that all vertices in $\{v_n\}$ are successful
in the process started from $(\zeta_1',\zeta_2')$. Hence
$\{v_n\}\cup \{e_k\}\subset \Psi^{\zeta_1',\zeta_2'}_1$, as
desired.$\hfill\Box$\medskip

\section{Proof of Theorem \ref{th:huvudresIII}}
\noindent We are now in a position to prove Theorem
\ref{th:huvudresIII}. The proof is based on a coupling of the
processes $P^{\lambda_1,\lambda_2}_{\xi_1,\xi_2}$ and
$P^{\lambda_1,\lambda_2}_{\xi_1',\xi_2'}$ that has certain
similarities with the coupling used in the proof of Proposition
1.1 in Deijfen \emph{et al} (2003). The geometrical aspects of our
proof are, however, quite different.\medskip

\noindent \emph{Proof of Theorem \ref{th:huvudresIII}:} Pick
finite sets $\xi_1,\xi_2,\xi_1',\xi_2'\subset\mathbb{Z}^d$ such
that the pairs $(\xi_1,\xi_2)$ and $(\xi_1',\xi_2')$ are fertile.
We will show that, if
$P^{\lambda_1,\lambda_2}_{\xi_1,\xi_2}(A)>0$, then
$P^{\lambda_1,\lambda_2}_{\xi_1',\xi_2'}(A)>0$ as well.
Interchanging the roles of $(\xi_1,\xi_2)$ and $(\xi_1',\xi_2')$
in the below arguments gives the reverse implication. As pointed
out before, by time-scaling and symmetry, it suffices to consider
the case when $\lambda_1=1$ and $\lambda_2\leq 1$, that is, when
the type 1 infection has rate 1 and is more powerful than the type
2 infection.\medskip

\noindent To begin with, we need some notation. For a finite set
$\eta\subset \mathbb{Z}^d$, let $\bar{m}(\eta)$ and
$\b{$m$}(\eta)$ be the points in $\mathbb{Z}^d$ with coordinates
$$
\bar{m}_i(\eta)=\max\{x_i;\hspace{0.1cm}x_i \textrm{ is the $i$th
coordinate of a point in }\eta\}
$$
and
$$
\b{$m$}_i(\eta)=\min\{x_i;\hspace{0.1cm}x_i \textrm{ is the $i$th
coordinate of a point in }\eta\}
$$
respectively, and define
$$
B_{\eta}=\{x\in\mathbb{Z}^d;\hspace{0.1cm}\b{$m$}_i(\eta)\leq
x_i\leq \bar{m}_i(\eta)\textrm{ for all }i=1,\ldots,d\},
$$
that is, $B_{\eta}$ is the smallest box that contains the set
$\eta$. Write $B_\eta^{+k}$ for the box $B_\eta$ enlarged by $k$
sites in each direction, that is,
$$
B_\eta^{+k}=\{x\in\mathbb{Z}^d;\hspace{0.1cm}\b{$m$}_i(\eta)-k\leq
x_i\leq \bar{m}_i(\eta)+k\textrm{ for all }i=1,\ldots,d\}.
$$
Finally, let $\xi=\xi_1\cup \xi_2\cup \xi_1'\cup \xi_2'$.\medskip

\noindent Now, first consider a process started from
$(\xi_1,\xi_2)$ and let $\tau$ be the time when the box
$B_\xi^{+2}$ is fully infected in this process, that is,
$$
\tau=\inf\{t;\hspace{0.1cm}B_\xi^{+2}\subset
\Gamma^{\xi_1,\xi_2}(t)\}.
$$
Then consider a process started from $(\xi_1',\xi_2')$, coupled
with the one started from $(\xi_1,\xi_2)$ in such a way that it
evolves independently up to time $\tau$ and then uses the same
Poisson processes as the process started from $(\xi_1,\xi_2)$ to
generate the infections after time $\tau$. The notation for this
coupled process is equipped with a hat-symbol, for instance
$\hat{\Gamma}^{\xi_1',\xi_2'}(t)$ denotes the set of infected
sites at time $t$. We will describe a scenario for the time
interval $[0,\tau]$ in the coupled process that -- combined with
Lemma \ref{lemma1} -- guarantees that both infection types grow
unboundedly in this process given that they do so in the process
started from $(\xi_1,\xi_2)$. To this end, assume that mutual
unbounded growth occurs in the process started from
$(\xi_1,\xi_2)$, that is, assume that both infection graphs
$\Psi_1^{\xi_1,\xi_2}$ and $\Psi_2^{\xi_1,\xi_2}$ contain an
infinite path. Let $x_1$ be the last site on $\partial
\Gamma^{\xi_1,\xi_2}(\tau)$ that is touched by an infinite path in
$\Psi_1^{\xi_1,\xi_2}$ and let $\tilde{x}_1$ be the last site on
$\partial B_\xi^{+2}$ that is touched by the path through $x_1$.
(Here ``last'' refers to time, that is, sites are ordered with
respect to the time when they are infected.) For simplicity we
assume that $\Gamma^{\xi_1,\xi_2}(\tau)$ has no holes so that
there are no uninfected sites completely surrounded by
$\Gamma^{\xi_1,\xi_2}(\tau)$. The desired scenario for the coupled
process is as follows:

\begin{itemize}
\item[1.] By definition of the box $B_{\xi_1'\cup \xi_2'}$, at time
zero there are at least two infected sites on $\partial
B_{\xi_1'\cup \xi_2'}$. If none of these is type 1 infected,
assume that some site on $\partial B_{\xi_1'\cup \xi_2'}$ becomes
type 1 infected before any type 2 infection takes place. This is
possible because of the assumption that $\xi_1'$ is not strangled
by $\xi_2'$.

\item[2.] Pick a type 1 infected point $x$ on $\partial B_{\xi_1'\cup
\xi_2'}$ and assume that the type 1 infection reaches $\partial
B_\xi^{+1}$ via the shortest possible path from $x$ without any
type 2 infections occurring. Also without any type 2 infections
occurring, suppose that the type 1 infection wanders the shortest
way along $\partial B_\xi^{+1}$ to the nearest neighbor of the
site $\tilde{x}_1\in\partial B_\xi^{+2}$ and then moves out to
$B_\xi^{+2}$ by infecting $\tilde{x}_1$. If $\tilde{x}_1$ happens
to be a corner point of $B_\xi^{+2}$, it does not have a nearest
neighbor on $B_\xi^{+1}$. To deal with this case, let
$\tilde{\tilde{x}}_1$ be the last (in time) non-corner point on
$B_\xi^{+2}$ that is touched by the infinite path through $x_1$.
Suppose then that the type 1 infection wanders along $\partial
B_\xi^{+1}$ to the nearest neighbor of $\tilde{\tilde{x}}_1$,
moves out to $B_\xi^{+2}$ by infecting $\tilde{\tilde{x}}_1$ and
then follows the infinite type 1 path to $\tilde{x}_1$.\medskip

\noindent The configuration at this stage of the construction
consists of the initial sets $(\xi_1',\xi_2')$ and a type 1 path
linking the set $\xi_1'$ to the point $\tilde{x}_1$ on $\partial
B_\xi^{+2}$; see Figure 1(a).

\begin{figure}
\centering \mbox{\subfigure[Configuration after Step 2. The dashed
line indicates the infinite type 1 path through $x_1$ in the
process started from $(\xi_1,\xi_2)$.
]{\epsfig{file=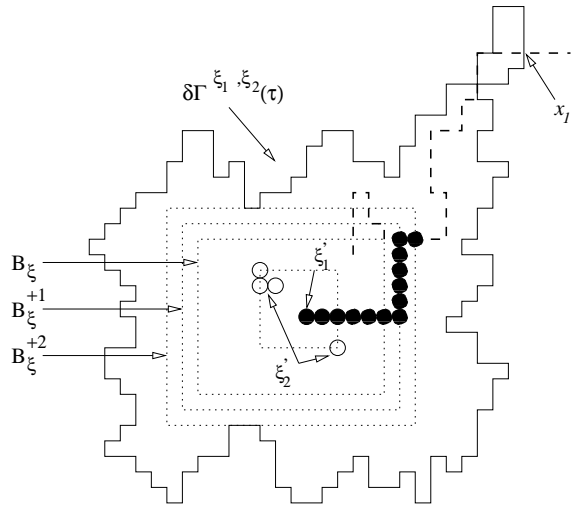,height=6.65cm}}}\par
\mbox{\subfigure[Configuration after Step 4.
]{\epsfig{file=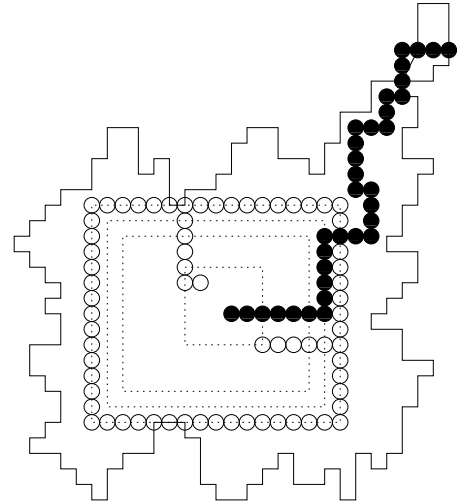,height=6.65cm}}}\caption{Development
of the infection in the process started from $(\xi_1',\xi_2')$.
Black circles represent type 1 infected sites and white circles
represent type 2 infected sites.}
\end{figure}

\item[3.] Now assume that the type 1 infection lies still while the
type 2 infection wanders out from $\xi_2'$ to $\partial
B_\xi^{+2}$ and invades all sites on $\partial B_\xi^{+2}$ except
$\tilde{x}_1$ (and possibly also $\tilde{\tilde{x}}_1$ and one or
more corner points on the type 1 path between $\tilde{x}_1$ and
$\tilde{\tilde{x}}_1$), which is already type 1 infected. Since
$\xi_2'$ is not strangled by $\xi_1'$, a path from $\xi_2'$ to
$\partial B_\xi^{+2}$ does indeed exist.

\item[4.] Suppose that the type 1 infection wanders from $\tilde{x}_1$
to the site $x_1$ on $\partial \Gamma^{\xi_1,\xi_2}(\tau)$ along
the infinite path through $x_1$ in $\Psi_1^{\xi_1,\xi_2}$ while no
type 2 infections occur.\medskip

\noindent To summarize, at this stage all sites on $\partial
B_\xi^{+2}$ except $\tilde{x}_1$ (and possibly also
$\tilde{\tilde{x}}_1$ and one or more corner points on the type 1
path between $\tilde{x}_1$ and $\tilde{\tilde{x}}_1$) are type 2
infected. We also have a type 1 path, passing $\partial
B_\xi^{+2}$ at $\tilde{x}_1$, that reaches from $\xi_1'$ to $x_1$;
see Figure 1(b).

\item[5.] Assume that $\partial \Gamma^{\xi_1,\xi_2}(\tau)$ is
infected as follows: If there are sites on $\partial
\Gamma^{\xi_1,\xi_2}(\tau)$ that cannot be reached from $\partial
B_\xi^{+2}$ using paths in $\Gamma^{\xi_1,\xi_2}(\tau)$ without
using sites on the type 1 path through $x_1$ -- this is the case
if $x_1$ is located on a cape as displayed in Figure 1 -- suppose
that the type 2 infection lies still until the type 1 infection
has invaded these sites, one at a time, starting from $x_1$. Then
assume that the rest of $\partial \Gamma^{\xi_1,\xi_2}(\tau)$ is
type 2 infected while no type 1 infections occur. Note that the
sites on $\partial \Gamma^{\xi_1,\xi_2}(\tau)$ that cannot be
reached without using sites on the type 1 path through $x_1$ must
constitute a connected subset of $\partial
\Gamma^{\xi_1,\xi_2}(\tau)$. Also, some thought reveals that these
sites must be type 1 infected in the process started from
$(\xi_1,\xi_2)$.

\item[6.] Suppose that the interior of
$\Gamma^{\xi_1,\xi_2}(\tau)$ is filled with infection without any
sites outside $\Gamma^{\xi_1,\xi_2}(\tau)$ being infected.

\item[7.] Let $T$ denote the time when the above scenario is
completed. Assume that $T\leq \tau$ and suppose that no infections
at all take place in the time interval $(T,\tau]$.
\end{itemize}

\noindent The above scenario clearly has positive probability
because it depends only on finitely many infections. Furthermore,
given this scenario, the following hold at time $\tau$:

\begin{itemize}
\item[--] $\hat{\Gamma}^{\xi_1',\xi_2'}(\tau)=\Gamma^{\xi_1,\xi_2}(\tau):=\Gamma$;
\item[--] $\Big(\Gamma_2^{\xi_1,\xi_2}(\tau) \cap
\partial\Gamma\Big)\subset\Big(\hat{\Gamma}_2^{\xi_1',\xi_2'}(\tau)\cap
\partial\Gamma\Big)$;
\item[--] $x_1\in\Gamma_1^{\xi_1,\xi_2}(\tau)\cap\hat{\Gamma}_1^{\xi_1',\xi_2'}
(\tau)\cap\partial\Gamma$.
\end{itemize}

\noindent After time $\tau$ the coupled process is based on the
same Poisson processes that were used to generate the process
started from $(\xi_1,\xi_2)$. It follows from Lemma \ref{lemma1}
with $(\zeta_1,\zeta_2)=
\big(\Gamma^{\xi_1,\xi_2}_1(\tau),\Gamma^{\xi_1,\xi_2}_2(\tau)\big)$
and $(\zeta_1',\zeta_2')=\big(\Gamma^{\xi_1',\xi_2'}_1(\tau),
\Gamma^{\xi_1',\xi_2'}_2(\tau)\big)$ that there is at least one
infinite path in both $\hat{\Psi}_1^{\xi_1',\xi_2'}$ and
$\hat{\Psi}_2^{\xi_1',\xi_2'}$. Hence we have mutual unbounded
growth in the coupled process.\medskip

\noindent Now let $A_{\xi_1,\xi_2}$ denote the event that both
infection types grow unboundedly in the process started from
$(\xi_1,\xi_2)$ and let $\hat{A}_{\xi_1',\xi_2'}$ denote the same
event in the coupled process. Trivially,

\begin{equation}\label{eq:4}
P(\hat{A}_{\xi_1',\xi_2'})\geq
P(\hat{A}_{\xi_1',\xi_2'}|A_{\xi_1,\xi_2})P(A_{\xi_1,\xi_2}).
\end{equation}

\noindent The above reasoning shows that, if both infection types
grow unboundedly in the process started from $(\xi_1,\xi_2)$, then
the scenario described in 1-7 guarantees that they do so in the
coupled process as well. Hence the first factor on the right-hand
side in (\ref{eq:4}) is positive. The last factor on the
right-hand side is positive by assumption. Thus
$P(\hat{A}_{\xi_1',\xi_2'})>0$ and, since
$P(\hat{A}_{\xi_1',\xi_2'})=P_{\xi_1',\xi_2'}^{1,\lambda_2}(A)$,
we are done.$\hfill\Box$\medskip

\section*{References}

\noindent Deijfen, M., H\"{a}ggstr\"{o}m, O. and Bagley, J.
(2003): A stochastic model for competing growth on $\mathbb{R}^d$,
\emph{Markov Proc. Relat. Fields}, to appear.\medskip

\noindent H\"{a}ggstr\"{o}m, O. and Pemantle, R. (1998): First
passage percolation and a model for competing spatial growth,
\emph{J. Appl. Prob.} \textbf{35}, 683-692.\medskip

\noindent H\"{a}ggstr\"{o}m, O. and Pemantle, R. (2000): Absence
of mutual unbounded growth for almost all parameter values in the
two-type Richardson model, \emph{Stoch. Proc. Appl.} \textbf{90},
207-222.\medskip

\noindent Newman, C.M. (1995): A surface view of first-passage
percolation, \emph{Proceedings of the 1994 International Congress
of Mathematicians}, ed. S.D. Chatterij, pp. 1017-1023,
Birkh\"{a}user.\medskip

\noindent Richardson, D. (1973): Random growth in a tessellation,
\emph{Proc. Cambridge Phil. Soc.} \textbf{74}, 515-528.\medskip

\end{document}